\documentclass[AMA,STIX1COL]{WileyNJD-v2}

\articletype{}

\usepackage{amsmath}
\usepackage{physics}

\usepackage{lineno}
\linenumbers

\let\oldequation\equation
\let\oldendequation\endequation

\renewenvironment{equation}
  {\linenomathNonumbers\oldequation}
  {\oldendequation\endlinenomath}

\let\oldalign\align
\let\oldendalign\endalign

\renewenvironment{align}
  {\linenomathNonumbers\oldalign}
  {\oldendalign\endlinenomath}

\newcommand{\tidx}{\mathrm{t}}
\newcommand{\fedsp}{\mathbf{u}}
\newcommand{\fevel}{\mathbf{\dot u}}
\newcommand{\feacc}{\mathbf{\ddot u}}
\newcommand{\cs}{c_\mathrm{s}}
\newcommand{\cl}{c_\mathrm{p}}
\begin{document}

\title{A three-dimensional hybrid finite element -- spectral boundary integral method for modeling earthquakes in complex unbounded domains}

\author[1,2]{Gabriele Albertini}
\author[3]{Ahmed Elbanna}
\author[1]{David S. Kammer}

\address[1]{\orgdiv{Institute for Building Materials}, \orgname{ETH Zurich},
\orgaddress{\country{Switzerland}}}

\address[2]{\orgdiv{School of Civil and Environmental Engineering}, \orgname{Cornell University}, \orgaddress{Ithaca, \state{NY}, \country{USA}}}

\address[3]{\orgdiv{Department of Civil and Environmental Engineering}, \orgname{University of Illinois}, \orgaddress{Urbana-Champaign, \state{IL}, \country{USA}}}

\corres{David S. Kammer \email{dkammer@ethz.ch}}

%
%

\abstract[Summary]{
We present a 3D hybrid method which combines the Finite Element Method (FEM) and the Spectral Boundary Integral method (SBIM) to model nonlinear problems in unbounded domains. The flexibility of FEM is used to model the complex, heterogeneous, and nonlinear part  -- such as the dynamic rupture along a fault with near fault plasticity -- and the high accuracy and computational efficiency of SBIM is used to simulate the exterior half spaces perfectly truncating all incident waves. The exact truncation allows us to greatly reduce the domain of spatial discretization compared to a traditional FEM approach, leading to considerable savings in computational cost and memory requirements. The coupling of FEM and SBIM is achieved by the exchange of traction and displacement boundary conditions at the computationally defined boundary. The method is suited to implementation on massively parallel computers. We validate the developed method by means of a benchmark problem. Three more complex examples with a low velocity fault zone, low velocity off-fault inclusion, and  interaction of multiple faults, respectively, demonstrate the capability of the hybrid scheme in solving problems of very large sizes. Finally, we discuss potential applications of the hybrid method for problems in geophysics and engineering.
}


%
%

\maketitle

\section{Introduction}\label{sec:intro}

Earthquakes are a prime example of complex natural processes with far-from-equilibrium nonlinear dynamics at multiple scales. 
The lack of quantitative data on timescales capturing multiple large earthquake cycles is a fundamental impediment for progress in the field. 
Physics-based simulations provide the only path for overcoming the lack of data and elucidating the multi-scale dynamics and spatio-temporal patterns that extend the knowledge beyond sporadic case studies and regional statistical laws.

The multiscale nature of the earthquake phenomena is manifested as follows. Spatially, a moderate-size earthquake typically propagates over tens of kilometres. 
However, the physical processes governing the rupture propagation operates within a narrow region at the rupture tip, called the process zone, which may not exceed a few millimetres in size if realistic laboratory-based friction parameters are used.\cite{noda_earthquake_2009} 
Temporally, an earthquake event, where rapid slip occurs, only lasts for few to tens of seconds.
However, the time span between successive large earthquakes may be tens to hundreds of years.\cite{lapusta_elastodynamic_2000} 
Thus, there exists approximately a decade of spatial and temporal scales that must be resolved in a target physics-based simulation of earthquakes and aseismic slip. 
This necessitates innovation in modeling both the fast dynamic rupture with extreme localization and the slow quasi-static slip, during the interseismic period, that exhibits gradual variations. 
This is a fundamental challenge in earthquake source physics which has been a focus of computational earthquake mechanics over the past four decades.

Historically, numerical methods for simulating earthquakes and aseismic slip may be classified broadly into two categories: boundary-based methods and domain-based methods. 
The boundary integral formulation enables reducing the spatial dimension of the problem by one, by invoking the representation theorem of linear elastodynamics, transforming 2D problems into 1D and 3D problems into 2D.\cite{cochard_dynamic_1994,geubelle_spectral_1995} 
The spectral formulation of the boundary integral equations has been transformative in seismic applications (e.g. Lapusta et al.\cite{lapusta_elastodynamic_2000} and references therein). 
For example, Lapusta et al.\cite{lapusta_elastodynamic_2000} derived accurate adaptive time-stepping algorithms and truncation of convolution integrals that enabled, for the first time, the consistent elastodynamic simulation of a long sequence of events combining rapid slip during earthquake ruptures and slow deformation during the interseismic periods.
Nonetheless, the method is limited to homogeneous linear elastic bulk.
While the method may be applied, in principle, to heterogeneous linear elastic materials, the lack of a closed form representation of the Green’s function either inhibits the method from providing a well-defined solution to many problems of interest or makes it less computationally attractive. 
Furthermore, the superior performance of the spectral approach and its computational efficiency is only possible for planar interfaces.
This precludes the representation of non-planar faults or direct incorporation of fault zone complexity (e.g. damage, and shear bands). 

\begin{figure}
    \centering
    \includegraphics[width=\textwidth]{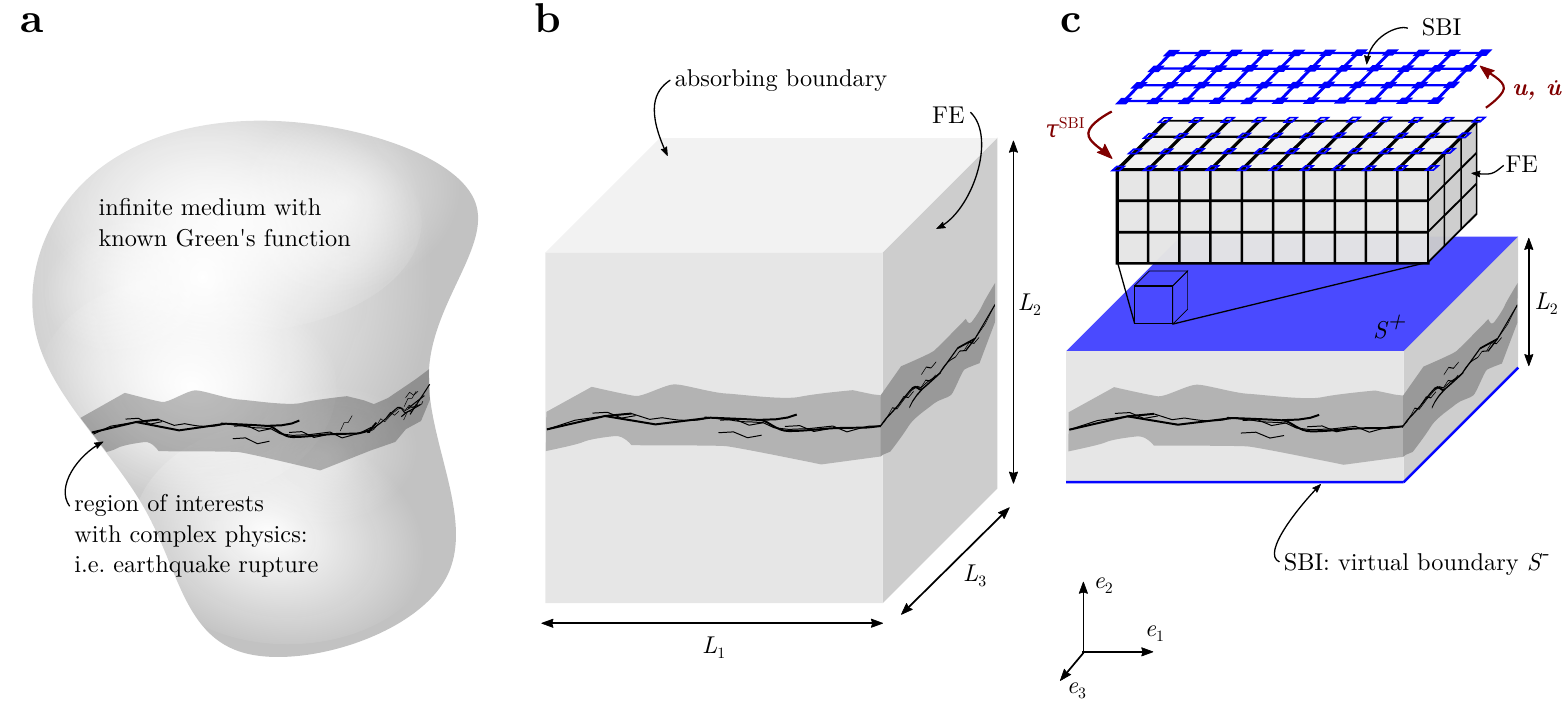}[h]
    \caption{Schematic representation of the physical problem (a), it's representation using a domain-based method such as FE (b), and using the hybrid method (c). The hybrid method couples a domain-based method with a boundary-based method, through the communication of nodal traction $\mathbf{\tau}^\mathrm{SBI}$,  displacement $\mathbf{u}$, and velocity $\mathbf{\dot{u}}$ at the boundaries of the virtual strip, $S^\pm$.}
    \label{fig:method}
\end{figure}

On the other hand, numerical methods based on bulk discretization such as the finite difference (FD) and finite element methods have been used in simulating earthquake ruptures since mid-1970s and early 1980s with the pioneering works of 
Boore et al.,\cite{boore_two-dimensional_1971}
Andrews,\cite{andrews_rupture_1976} 
Das \& Aki,\cite{das_numerical_1977} 
Archuleta \& Day,\cite{archuleta_dynamic_1980} 
Day, \cite{day_three-dimensional_1982} 
Virieux \& Madariaga, \cite{virieux_dynamic_1982} and others. 
These methods are more flexible than the boundary integral approaches in handling heterogeneities, nonlinearities, and fault geometry complexities (see Fig.~\ref{fig:method}a\&b). 
In recent years, highly accurate formulations were introduced, including the spectral finite element,\cite{komatitsch_introduction_1999,ampuero_etude_2002,festa_influence_2006,ma_radiated_2006,kaneko_spectral_2008}  
the discontinuous Galerkin method,\cite{kaser_arbitrary_2006,benjemaa_dynamic_2007,puente_dynamic_2009,pelties_three-dimensional_2012,tago_3d_2012}  
and higher-order FD schemes.\cite{cruz-atienza_3d_2007,dalguer_staggered-grid_2007,kozdon_simulation_2013} 
A main computational challenge of these methods is the need to discretize the whole bulk, which increases the computational demand by at least one order of magnitude compared to the boundary integral formulation. 
Furthermore, the computational domain must be truncated at a sufficient distance from the fault surface such that it would not affect the physical solution. 
While domain truncation has been achieved by the introduction of several widely-used absorbing boundary conditions such as boundary viscous damping,\cite{lysmer_finite_1969} perfectly matching layers,\cite{berenger_perfectly_1994} and infinite elements,\cite{bettess_infinite_1977} these methods have limitations.
Specifically, in all these methods, artificial reflections exist to varying degrees and the absorbing surfaces must be taken sufficiently far from the fault surface to ensure solution accuracy. 
Moreover, attempts to perform cycle simulations using these volume-based methods are rare and have been restricted mainly to the quasi-dynamic limit.\cite{erickson_bimaterial_2016} 
This is partially due to the high spatial discretization cost and the lack of a systematic approach to handle both dynamic and quasi-dynamic calculations in the same framework which is required for simulating both earthquake ruptures and intersesismic slow deformations. 
Another challenge in these methods is defining fault loading. 
Currently, this is done by applying displacement-controlled loading at the far boundaries of the simulation box. 
This, however, makes the fault stressing rate dependent on where the domain is truncated.
This problem is solved approximately in the SBI formulation by loading the fault directly through back-slip.

Both bulk and boundary approaches have their merits and limitations. 
The limitations are evident in 3D simulations where computational complexity grows like the element size to the fourth power rendering high resolution models a computational bottleneck. 
To that end, this paper proposes a new hybrid numerical scheme, for the full three dimensional elastodynamic problem, that combines the 3D FE method and the 2D SBI equation method to efficiently model fault zone nonlinearities and heterogeneities with high resolution while capturing large-scale elastodynamic interactions in the bulk. 
The main idea of the method is to enclose the heterogeneities in a virtual strip that is introduced for computational purposes only (see Fig.~\ref{fig:method}c). 
This strip is discretized using a volume-based numerical method, chosen here to be the finite element method due to its popularity and flexibility in handling complex geometry and arbitrary bulk heterogeneities. 
The top and the bottom boundaries of the virtual strip are handled using the independent SBI formulation\cite{geubelle_spectral_1995} with matching discretization. 
The coupling between the two methods is achieved through enforcing continuity of displacement and traction at the virtual boundaries.
The current work extends recent work by the authors and their groups over the past few years which first developed the hybrid scheme for the 2D dynamic anti-plane problem combining finite difference and spectral boundary integral methods, \cite{hajarolasvadi_new_2017} and the 2D dynamic in-plane problem using the finite element method for bulk discretization in the hybrid scheme.\cite{ma_hybrid_2019} 
Prior work has demonstrated the accuracy and computational efficiency of the coupled approach and its potential for modeling dynamic ruptures with high resolution fault zone physics\cite{ma_dynamic_2019} as well as extension to the quasi-dynamic limit and cycle simulations. \cite{abdelmeguid_novel_2019} 
The current extension to the full three dimensional case represents the culmination of these efforts. 

The remainder of this paper is organized as follows. In Section~\ref{sec:method}, we describe the physical model (Section~\ref{sec:physics}), and the numerical methods to solve it, which includes the finite-element method (Section~\ref{sec:fem}), the spectral boundary integral method (Section~\ref{sec:sbi}), and their coupling -- the hybrid method (Section~\ref{sec:hybrid}).
In Section~\ref{sec:TPV3}, we validate the hybrid method using the benchmark problem TPV3\cite{day_comparison_2005} of the Southern California Earthquake Center.
Next, we demonstrate the capabilities of the new hybrid method on more complex problems. We consider a low velocity fault zone in Section~\ref{sec:LVFZ},
a low velocity inclusion at a distance from the fault in Section~\ref{sec:HVFZ}, and interacting faults in Section~\ref{sec:step-over}.
Finally, we discuss the advantages of the hybrid method in terms of computational cost in Section~\ref{sec:discussion} and draw conclusions in Section~\ref{sec:conclusion}.

\section{Method}\label{sec:method}

\subsection{Physical Model}\label{sec:physics}

We solve the fully dynamic three-dimensional problem of a rupture propagating along a fault embedded in an elastic solid.
The conservation of linear momentum within the elastic domain $\Omega$ is given by
\begin{equation}
  \rho  \ddot{u}_i - \pdv{\sigma_{ij}}{x_j} = 0 \quad \text{in} \; \Omega
\label{eq:consmoment}
\end{equation}
where $\rho$ is the material density, $u_i$ the displacement vector, with the ``dot'' being the derivative with respect to time $t$, $\sigma_{ij}$ the Cauchy stress tensor, and $x_j$ the coordinate vector. Body forces are neglected.
Dirichlet boundary conditions are applied on  $S_u$ and Neumann boundary conditions are applied on $S_T$
\begin{align}
    u_i&=\bar{u}_i \qquad \text{on } S_u\label{eq:dirichletBC}\\\
    \sigma_{ij}n_j&=\bar \tau_i \qquad \text{on } S_\tau,\label{eq:neumannBC}
\end{align}
where $n_i$ is the normal vector to the surface $S_\tau$. Initially, the domain is assumed to be in equilibrium, and, hence, the initial conditions are given by $u_i(0) = u_i^0$ and $\dot{u}_i(0) = 0$.
We assume linear elastic material behavior:
\begin{equation}
  \sigma_{ij} = \lambda \, \delta_{ij} \varepsilon_{kk} + 2 \mu \, \varepsilon_{ij}
\label{eq:linelasticity}
\end{equation}
with the infinitesimal strain tensor $\varepsilon_{ij} = (\partial u_i / \partial x_j + \partial u_j / \partial x_i) / 2$  and the Lam\'e parameters $\lambda$, $\mu$ describing the elastic properties of the material.

The fault transmits stresses from one half-space to the other through interface tractions. 
We focus on tangential (friction) interaction and impose non-penetration/non-opening conditions to the normal component of the fault surfaces $S_f^\pm$. The local slip vector, which corresponds to the tangential fault opening vector, is given by 
\begin{equation}
\delta_i = R_{ij} (u_j^+ - u_j^-) \quad \text{on } S_f^\pm\,,    
\end{equation}
where $R_{ij}$ is the global-to-local rotation matrix and $+$ and $-$ indicate the upper and lower fault sides, respectively. Local slip $\delta$ is the amplitude of $\delta_i$. The fault is governed by a stick-slip behavior that is described by two states. A sticking section of the fault is described by:
\begin{equation}
  \dot \delta = 0 \quad \text{and} \quad \tau \leq \tau^s
  \label{eq:stick}
\end{equation}
where $\dot \delta$ is slip rate, $\tau$ the amplitude of the fault shear traction vector and $\tau^s$ the fault strength. A sliding fault section is described by:
\begin{equation}
  \dot \delta > 0 \quad \mathrm{and} \quad \tau = \tau^s ~.
  \label{eq:slip}
\end{equation}
A constitutive law is applied to model the fault strength evolution. For simplicity, we apply a linear slip-weakening friction law, \cite{ida_cohesive_1972} which is given by:
\begin{equation}
\mu (\delta) = 
\left \{ \begin{array}{ll}
  \mu_s - (\mu_s - \mu_k) \, \delta/\delta_c & \quad \text{for} \;  \delta < \delta_c\\
  \mu_k & \quad \text{for} \; \delta \geq \delta_c
\end{array}
\right .
\label{eq:slipweak}
\end{equation}
where $\mu_s$ and $\mu_k$ are the static and kinetic friction coefficient, respectively, and $\delta_c$ is the characteristic slip length to reach residual strength. The fault strength is then given by $\tau^s=\sigma^n \mu(\delta)$, where $\sigma^n$ is the normal stress.
Other friction laws, such as rate-and-state friction, \cite{dieterich_modeling_1979,ruina_slip_1983}
could also be applied in a similar framework, as shown by Kaneko et al.\cite{kaneko_spectral_2008}

\subsection{Finite Element Method (FEM)}\label{sec:fem}
The finite element method is based on a variational formulation of the governing equation and applies a discretization based on shape functions to find an approximate solution to the physical problem presented in Sec.~\ref{sec:physics}. A detailed description can be found in standard textbooks.\cite{belytschko_nonlinear_2013} 
The FEM approach transforms the strong form, \emph{i.e.,} the governing equation (\ref{eq:consmoment}),  to the  weak form by multiplying it with the test functions $\hat u(x)$, integrating it over the domain, and applying Green's identity, which results in:
\begin{equation}
    \int_\Omega \rho \ddot{u}_i \, \hat u_i \, \dd \Omega
    +\int_\Omega \sigma_{ij}\pdv{\,\hat u_i}{x_j} \dd \Omega
    - \int_{S_\tau}     \bar \tau_i \, \hat u_i \, \dd S 
    - \int_{S_f}     \tau_i \, \hat u_i \, \dd S 
    =0\,.
\end{equation}
Test functions are chosen smooth enough, such that all steps are well defined and vanish on the Dirichlet boundary. 
By choosing suitable interpolation functions, $N_I(x_J)=\delta_{IJ}$, and test function $\hat u(x)=\sum_I N_I(x)u_I(t)$, where the subscript $I$ represent the node index, $u_I(t)$ becomes the nodal displacements.
Using this standard FE approach, the weak form can be expressed as the following matrix equation
\begin{equation}
  \mathbf{M} \feacc + \mathbf{K} \fedsp  - \mathbf{f} - \mathbf{B \tau} = \mathbf{0}
\end{equation}
where $\mathbf{\ddot u}$ denotes the second time derivative of the displacement vector, $\mathbf{M}$ and $\mathbf{K}$ are the mass and stiffness matrix, respectively, $\mathbf{B}$ is a fault rotation-area matrix, $\mathbf{\tau}$ is the fault traction vector, and $\mathbf{f}$ is the force vector from Neumann boundary conditions.

We apply an explicit central-difference time integration formulation with a predictor-corrector formulation. The step-by-step procedure follows:
\begin{align} 
    \fevel^{pred}_{\tidx+1} &= \fevel_{\tidx} + \Delta t ~ \feacc_{\tidx} \label{eq:fepred}\\
    \fedsp_{\tidx+1} &= \fedsp_{\tidx} + \Delta t ~ \fevel^{pred}_{\tidx+1} \label{eq:feintdisp}\\
    \Delta \feacc &=\left(-\mathbf{Ku}_{\tidx +1} + \mathbf{f}  +\mathbf{B\tau}_{\tidx+1} \right)\mathbf{M}^{-1} -\feacc_{\tidx} \label{eq:tintf}\\
    \fevel_{\tidx+1}&=\fevel^{pred}_{\tidx+1} + \frac{1}{2} \Delta t \, \Delta \feacc\label{eq:fecorr}\\
    t&=t+\Delta t
\end{align}
where the subscript indicates the time step and $\Delta t$ is the current incremental time step, which is required to satisfy the Courant-Friedrichs-Lewy condition.\cite{courant_uber_1928} We apply a lumped mass matrix, which simplifies the computation of the inverse mass matrix and reduces computational cost of the time-integration scheme. The fault rotation matrix is scaled by the fault-surface area associated with each fault-split-node and thus transforms the fault traction vector $\mathbf{\tau}$ to a nodal force vector.

The fault traction vector $\mathbf{\tau}_{\tidx+1}$ in Eq.~(\ref{eq:tintf}) is computed by a forward Lagrange multiplier method,\cite{carpenter_lagrange_1991} which uses a prediction procedure to pre-computes the slip rate for the next time step. A similar approach was applied in spectral-element simulations. \cite{kaneko_spectral_2008}
We denote a fault discontinuity as $[|A|] = (A_+ - A_-)$ where subscript $+$ and $-$ indicate the upper and lower fault sides, respectively. The predicted slip rate during the next time step if no fault tractions were applied is given by
\begin{equation}
    [|\fevel_{\tidx+3/2}|] = [| \fevel^{pred}_{\tidx+1} - \frac{\Delta t}{2} \feacc_{\tidx} - \Delta t ~ \mathbf{M}^{-1} \left(\mathbf{K} \fedsp_{\tidx+1} - \mathbf{f}\right) |] ~.
  \label{eq:freeslip}
\end{equation}
We use the slip-rate predictor $[|\fevel_{\tidx+3/2}|]$ because the no-slip-rate condition $[| \fevel_{\tidx+3/2}|] = 0$ will ensure that the interface remains stuck and hence $[|\fedsp_{\tidx+2}|] = [|\fedsp_{\tidx+1}|]$. Using Eq.~\ref{eq:freeslip}, we can compute the traction required to maintain slip and impose stick condition on the fault by
\begin{equation}
  \mathbf{\tilde \tau}_{\tidx+1} = \frac{1}{2} \mathbf{Z} [|\fevel_{\tidx+3/2}|] ~,
  \label{eq:sticktraction}
\end{equation}
where $\mathbf{Z}$ is the fault impedance matrix given by $\mathbf{Z}^{-1} = \Delta t \left( \mathbf{M}_+^{-1} \mathbf{B}_+ + \mathbf{M}_-^{-1} \mathbf{B}_- \right)/ 2 $ and the following fault traction balance was applied $\mathbf{\tau} = - \mathbf{\tau}_+ = \mathbf{\tau}_-$.
The actual fault traction is computed by applying the stick-slip conditions given by Eq.~(\ref{eq:stick}) and Eq.~(\ref{eq:slip}): 
\begin{equation}
  \tau_{\tidx+1} = 
\left \{ \begin{array}{llll}
  \tilde \tau_{\tidx+1} & \mathrm{if} \quad \tilde \tau_{\tidx+1} \leq \tau^s_{\tidx+1} &\mathrm{(stick)}\\
  \tau^s_{\tidx+1} & \mathrm{otherwise} &\mathrm{(slip)}
\end{array}
\right .
\label{eq:stickslip}
\end{equation}
where $\tau_{\tidx+1}$ and $\tilde \tau_{\tidx+1}$ are individual entries in $\mathbf{\tau}_{\tidx+1}$ and $\mathbf{\tilde \tau}_{\tidx+1}$, respectively, and $\tau^s_{\tidx+1}$ is the fault strength at each split-node (node indicator is omitted for simplicity) and is governed by Eq.~(\ref{eq:slipweak}).

\subsection{Spectral Boundary Integral Method (SBIM)}\label{sec:sbi}

Boundary integral methods have the advantage of modeling the wave propagation problem in the entire domain $\Omega$ by using an integral relationship (in space and time) between the displacements and the tractions along the boundary of the domain $\partial \Omega$.
The advantage lies  in reduced computational cost and increased accuracy with respect to a finite-element or finite-difference method. 
For these reasons, boundary integral methods have been used extensively since the mid-1980s to study crack propagation problems.\cite{das_numerical_1980,andrews_dynamic_1985,boatwright_seismic_1986,das_numerical_1987,das_investigation_1988,israil_two-dimensional_1990,koller_modelling_1992,liu_hypersingular_1993,bonnet_regularization_1993,cochard_dynamic_1994,andrews_dynamic_1994}  

Consider the displacements and tractions at a the boundary of an semi-infinite half space with the boundary lying on the $e_1,e_3$ plane and the domain being infinite in the $e_2$ direction. 
Following the process described in Geubelle and Rice,\cite{geubelle_spectral_1995} the elastodynamic response of a 3D elastic half space is given by  
\begin{equation}\label{eq:SBI}
    \tau_i^{\mathrm{SBI}}(x_1,x_3,t) = \tau_i^\infty(x_1,x_3,t) - \eta_{ij}\frac{\mu}{\cs} \dot u_j(x_1,x_3,t) + s_i(x_1,x_3,t) ~,
\end{equation}
where $\eta_{ij}$ is a diagonal matrix with $\eta_{11}=\eta_{33}=1$ and $\eta_{22}=\cs/\cl$. $\cl$ and $\cs$ are the longitudinal and shear wave speeds of the material, respectively.
Eq.~(\ref{eq:SBI}) states that the traction at the surface of the half space, $\tau^\mathrm{SBI}_i$, equals the far field traction, $\tau^\infty_i$, plus a ``radiation damping'' term, $\eta_{ij}\frac{\mu}{\cs}\dot u_j^\pm$, and a spatiotemporal integral term $s_i$.
In this formulation the elastodynamic response is separated between local and nonlocal contributions.
$s_i$ represents the nonlocal elastodynamic long-range interaction between different parts of the surface, and the local effect $\eta_{ij}\frac{\mu}{\cs}\dot u_j$ represents wave radiation from the surface.

We use the spectral approach\cite{breitenfeld_numerical_1998} for computing $s_i$ which involves a Fourier transform in space and a convolutions in time, where the displacement history is convolved with the elastodynamic kernels. Please refer to independent formulation in \cite{breitenfeld_numerical_1998} for the derivation of the kernels and details for computing the nonlocal term, $s_i$.

\subsection{Hybrid FEM-SBIM Method}\label{sec:hybrid}

The hybrid method consists in coupling the FEM and the SBIM at the boundaries $S^\pm$, where the FE-domain is truncated (see Fig.~\ref{fig:method}c). At $S^\pm$, which we refer to as the virtual boundary, we apply an exact elastodynamic transparent boundary condition using the SBIM, which accounts for wave propagation in the infinite half-space beyond the FE truncation.
Depending on the FE scheme, Neumann or Dirichlet boundary conditions might be more suitable. For example, Dirichlet boundary conditions might result in a more stable algorithm.
Here, we present the Neumann approach (for the Dirichlet approach please refer to Ma et al. \cite{ma_hybrid_2019}).
We impose continuity condition at the boundaries $S^\pm$, which results in the FE force $\mathbf{f}$ being equal to the SBI traction $\tau^{\mathrm{SBI}}$ multiplied by a rotation-area matrix $\mathbf{B}^{\mathbf{SBI}}$.

The Neumann approach consists in solving the boundary integral relation, Eq.~(\ref{eq:SBI}), by using the displacements and velocity computed from the FEM. The resulting traction is then applied as a Neumann boundary condition in the FEM.
A time step of the hybrid method is computed as follows:
\begin{enumerate}
    \item FE compute explicit time integration Eq.~(\ref{eq:feintdisp}) and  predict velocity Eq.~(\ref{eq:fepred})
    
    \item copy $\fedsp_{\tidx+1}$ and $\fevel^{pred}_{\tidx+1}$ from FE to SBI
    
    \item SBI compute response of half space, $\tau^\mathrm{SBI}$, for given displacement history $u(x_1,x_3,t)$ and current velocity prediction $\dot{u}^{pred}_{\tidx+1}$ Eq.~(\ref{eq:SBI})
    
    \item apply SBI interface traction as Neumann boundary condition in FE: $\mathbf{f=B^{\mathrm{SBI}}} \mathbf{\tau}^\mathrm{SBI}$ 
    
    \item FE compute friction traction $\mathbf{\tau_{\tidx+1}}$ using Eq.~(\ref{eq:stickslip})
    
    \item FE compute acceleration increment Eq.~(\ref{eq:tintf})

    \item FE correct velocity Eq.~(\ref{eq:fecorr})

\end{enumerate}

Alternative coupling methods, \textit{e.g.}, Langrange multiplier, could also be applied. However, as we will show in Sec.~\ref{sec:TPV3}, the simple staggered approach proposed here provides excellent accuracy and is optimal in terms of computational efficiency. 
Further, in the finite-element domain, we apply 8-node linear hexagonal elements in a regular mesh in all presented problems.

\section{Benchmark problem TPV3: Earthquake rupture in unbounded homogeneous domain}\label{sec:TPV3}
\subsection{Setup}

\begin{figure}[h]
    \centering
    \includegraphics[width=\textwidth]{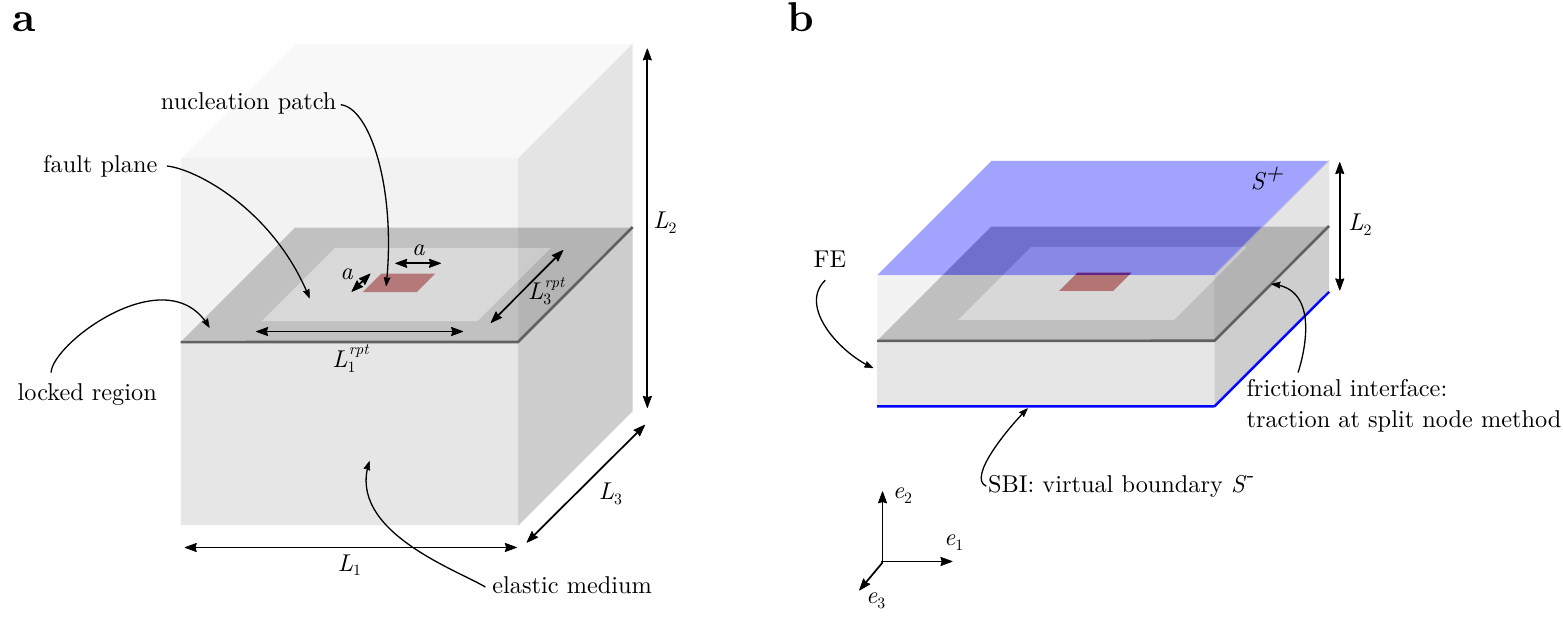}
    \caption{Setup of benchmark problem TPV3. 
    (a)~Earthquake rupture in unbounded domain with nucleation over square region of size $a=3~\text{km}$ and fault regions of size $L_1^{rpt}=30~\text{km}$ and $L_3^{rpt}=15~\text{km}$.
    (b)~Hybrid setup: FE domain with SBI as elastodynamic boundary condition.}
    \label{fig:TPV3_setup}
\end{figure}

We verify the hybrid method with the benchmark problem TPV3 from the SCEC Dynamic Rupture Validation exercises (\url{https://strike.scec.org/cvws/}).
The problem considers a planar fault, governed by linear slip-weakening friction, embedded in a homogeneous linear elastic bulk (see Fig.~\ref{fig:TPV3_setup}a). 
The elastic bulk has a density of $\rho=2670~\text{kg/m}^3$,
pressure wave speed $c_p=6000~\text{m/s}$,
and shear wave speed $c_s=3464~\text{m/s}$.
The friction properties are uniform and characterized by $\mu_s=0.677$, $\mu_k=0.525$ and $d_c=0.4~\text{m}$.
A uniform background shear, $\tau_0=70~\text{MPa}$, and normal stress, $\sigma_0=120~\text{MPa}$, are applied. The rupture is nucleated at the center of the fault over a square patch of size $a^2$ by instantaneously increasing the shear stress to a value higher than the static friction.
After nucleation, the rupture quickly propagates across the entire fault.

Note that this problem does not present any off-fault non-linearities and the fault is planar. Hence, the hybrid method is not required for this particular problem, which could be solved solely by the SBIM. However, we use this problem to verify the hybrid method by comparing the results with the reference solution of the SBIM. Additional examples, which do include non-linearities that require the hybrid method, will be presented in the following sections. 

Additionally, in this example (as well as in Sec.~\ref{sec:LVFZ} and \ref{sec:HVFZ}) the fault  (\emph{i.e.}, the 13 plane) represented a symmetry plane. The implication of this are twofold: (i) the normal stress on the fault remains constant throughout the simulation and so do the peak and residual friction strengths (ii) there is no need of explicitly modelling the bottom half space of the simulation domain (\emph{i.e.,} $x_2<0)$), which is taken into account for when applying the frictional traction on the fault by enforcing $[|u_2|]=0$ and $\sigma_{22}=const$.
Note that, by applying the hybrid method, we model the entire top (\emph{i.e.,} $x_2>0$) half space. However,  there are two additional symmetry planes: the 12 plane and the 23 plane. Thus, the FE computational domain could be further reduced by using the appropriate boundary conditions.

\subsection{Results}

\begin{figure}[h]
    \centering
    \includegraphics[width=\textwidth]{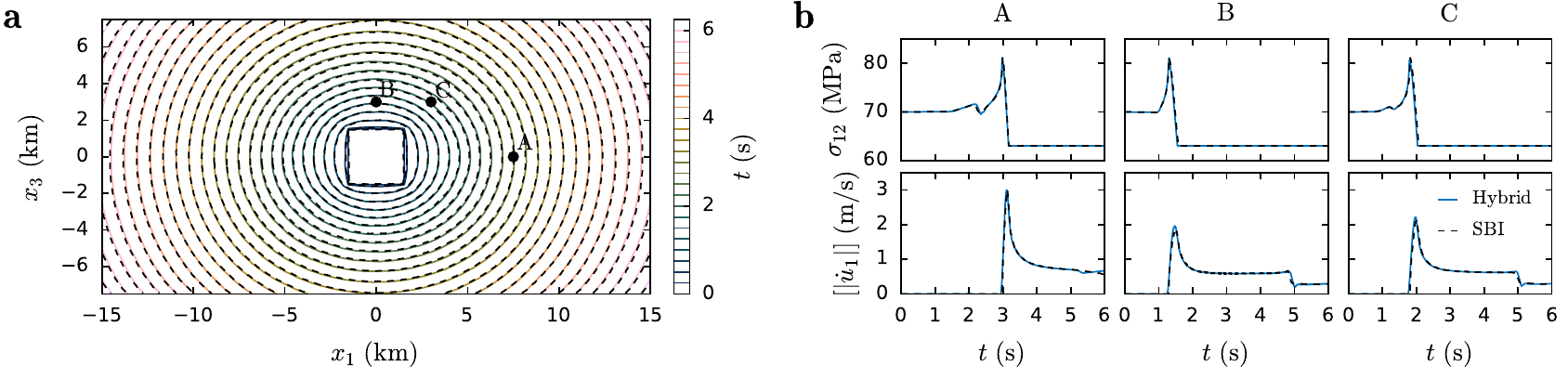}
    \caption{Benchmark problem TPV3 solved using hybrid method  (blue lines) and using the SBIM (dashed black lines).
    (a)~Contour of rupture front position each 0.5s.
    (b)~Fault shear stress, $\sigma_{12}$, and slip rate, $\dot{u}_1$, at three stations A, B, and C with position shown in (a).}
    \label{fig:TPV3_results}
\end{figure}

\begin{figure}[h]
    \centering
    \includegraphics[width=0.75\textwidth]{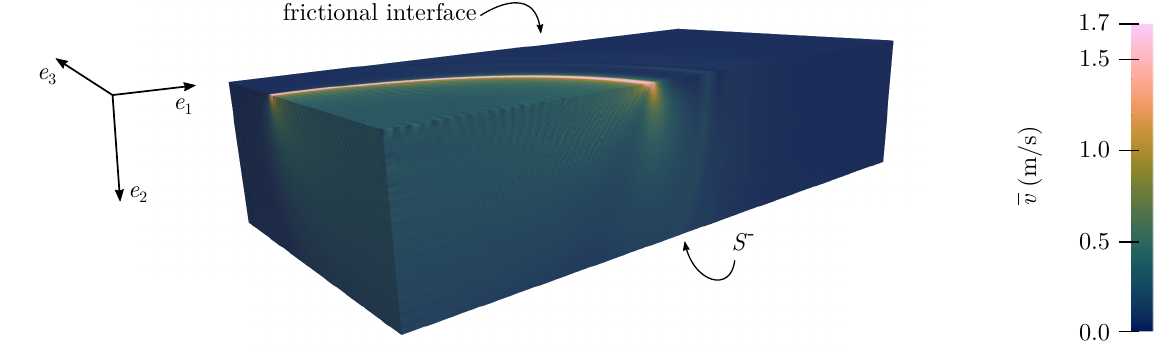}
    \caption{Benchmark problem TPV3 solved using hybrid method. Velocity magnitude field, $\bar v$, at $t=3.4~\text{s}$. For better visualization we applied a much thicker virtual strip, \textit{i.e.}, $L_2=10~\text{km}$ instead of $L_2=0.2~\text{km}$ as applied for simulations shown in Fig.~\ref{fig:TPV3_results}.}
    \label{fig:TPV3_para}
\end{figure}

We present the results of the hybrid method with a virtual strip width $L_2=4\Delta x$ (see Fig.~\ref{fig:TPV3_setup}b) and compare it with the reference solution. For both methods we use the same spatial discretization with $\Delta x= 50~\text{m}$.
Rupture front position (see Fig.~\ref{fig:TPV3_results}a), and stress and slip time history at three stations (see Fig.\ref{fig:TPV3_results}b) show excellent agreement between the hybrid method and the SBIM.
Fig.~\ref{fig:TPV3_para} shows the shear stress $\sigma_{12}$ field on one quadrant of the virtual strip. 
The rupture front is characterized by an abrupt change from peak to residual strength. 
The excellent agreement with the reference solution (Fig.~\ref{fig:TPV3_results}) demonstrates that the the elastodynamic boundary condition enforced on the planes $S^\pm$ does not cause any artificial wave reflection even though the virtual strip is extremely thin.

\begin{figure}[h]
    \centering
    \includegraphics[width=\textwidth]{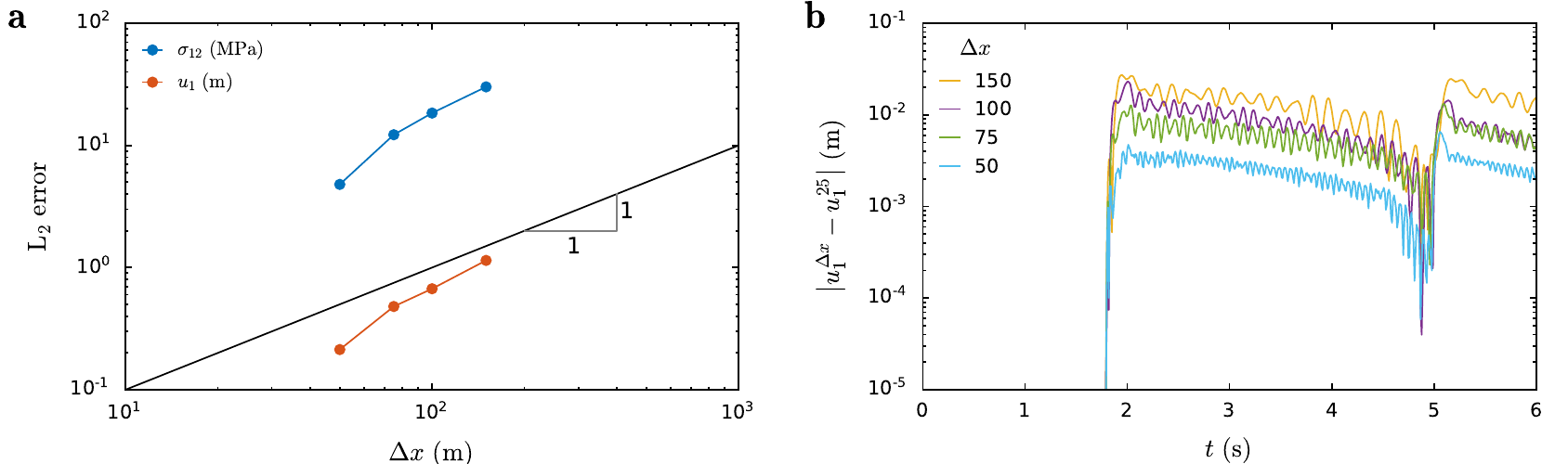}
    \caption{Benchmark problem TPV3 solved with hybrid method -- mesh convergence study 
    (a)~$\mathrm{L_2}$ error as function of mesh size $\Delta x$ computed over the entire fault at $t=3~\text{s}$. 
    (b)~Absolute error as function of time at station C.}
    \label{fig:TPV3_error}
\end{figure}

We perform a mesh convergence study of the hybrid method and show that the L$_2$ norm of the error, computed using a reference solution with $\Delta x=25~\text{m}$,  decreases linearly (see Fig.~\ref{fig:TPV3_error}a). 
The hybrid method combines linear finite elements with a higher precision spectral boundary method. Therefore, the convergence rate of the hybrid method corresponds to the rate of the least accurate of the methods it combines, \emph{i.e.,} the linear finite element method. Hence, the hybrid method does not loose any accuracy compared to a fully FEM model.

Since we are dealing with a dynamic problem, we also show the temporal evolution of the error at station C (see Fig.~\ref{fig:TPV3_error}b). The error is initially zero because the waves and the rupture have not reached the station yet. At $t\approx 1.8~\text{s}$, the rupture reaches the station (see also Fig.~\ref{fig:TPV3_results}b) and hence the error increases rapidly. It then remains approximately constant while the fault continues to slide until reflected waves from the boundary between rupture region and locked region (not the virtual boundary) reach the station. At this point, we notice a temporary drop in the error before it increases again to the same level of error observed before. Overall, we find that the error remains mostly constant over the duration of the simulation and decreases with mesh refinement.

Even though this benchmark problem is linear elastic and does not necessitate the use of the hybrid method, it illustrates its capability of efficiently and accurately truncating elastic waves in the vicinity of the fault with no artificial reflections from the virtual boundaries, $S^\pm$, which were only two elements away from the fault. 
In more complex scenarios, this virtual strip might need to be larger in order to fully describe the source of non-linearities or heterogeneity. Nevertheless, this efficient near-field truncation algorithm enables us to decrease the domain of finite-element discretization, compared to a fully FEM model, and apply a volumetric mesh only in a narrow strip around the fault, which results in considerable savings in terms of both computational time and memory, as we will discuss further in Sec.~\ref{sec:discussion}.

\section{Earthquake rupture with LVFZ: pulse-like behaviour}\label{sec:LVFZ}

\subsection{Setup}

The previous example was a benchmark problem and could have been solved by a boundary-element approach without any discretization of the bulk. Hence, the hybrid method was not required. In the following, we will consider more complex problems, which require volumetric discretization. 
First, we consider a slip-weakening fault with a low velocity fault zone (LVFZ).
LVFZ are found in most mature faults, where the near fault rock is considerably damaged and, as a consequence, has a reduced wave speed ranging from 20\% to 60\% with respect to the host rock.\cite{ma_effect_2015,huang_pulse-like_2011,huang_earthquake_2014,huang_potential_2016,albertini_off-fault_2017}
In 2D setups, when the reduction is high enough, the rupture behaves like a pulse. The results presented here will confirm this behavior on a 3D setup.

\begin{figure}[h]
    \centering
    \includegraphics[width=\textwidth]{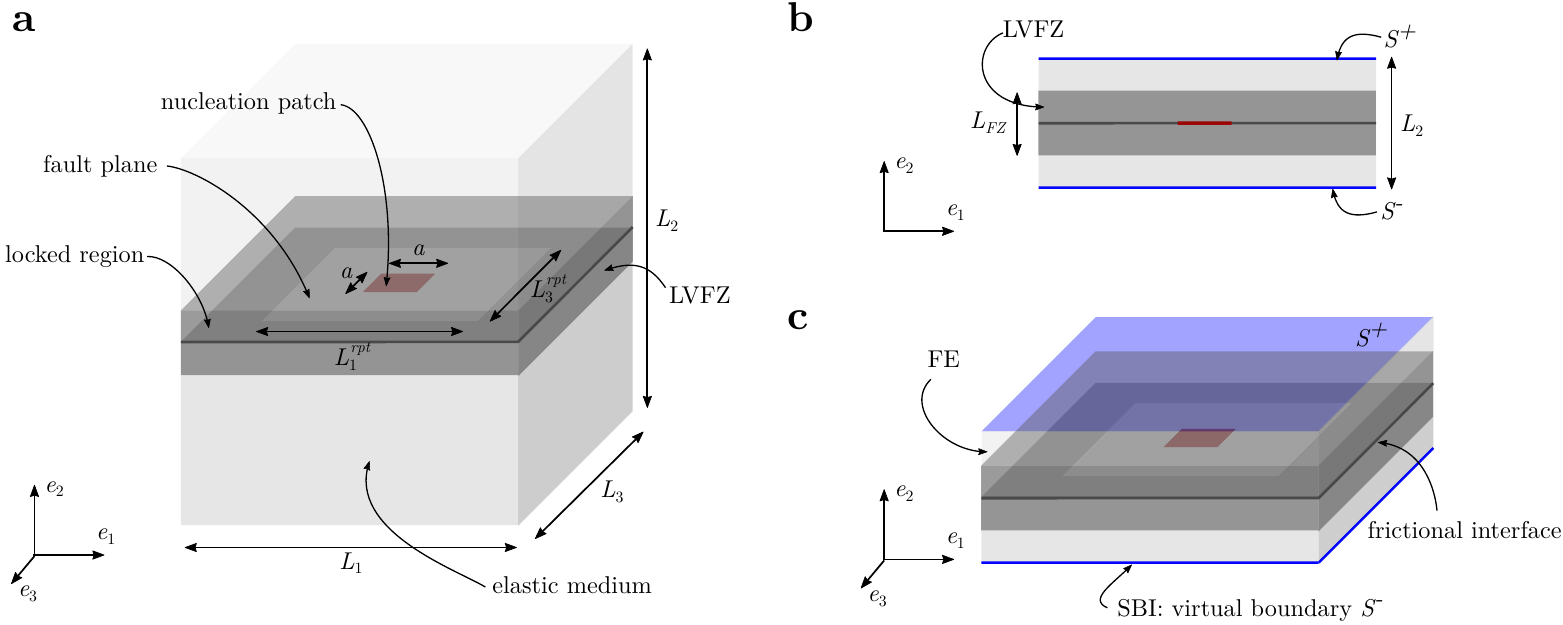}
    \caption{Setup of example earthquake rupture in unbounded domain with LVFZ.
      (a)~Fault plane geometry and nucleation patch are analogous to the previous example but with a larger size: $L_1^{rpt}=60~\text{km}$, $L_3^{rpt}=30~\text{km}$, and $a=3.2~\text{km}$.
      The fault zone region surrounding the fault is more compliant and presents a thickness of $L_{FZ}=1.6~\text{km}$.
    (b)~Hybrid setup: FE virtual strip with $L_2=2~\text{km}$ and with SBI as elastodynamic boundary condition.}
    \label{fig:LVFZ_setup}
\end{figure}

We consider a velocity reduction of 20\% with respect to the surrounding host rock, which has the same elastic properties as in Sec.~\ref{sec:TPV3}. The fault geometry is given in Fig.~\ref{fig:LVFZ_setup}, and is governed by linear slip-weakening friction with $\mu_s=0.677$, 
$\mu_k=0.564$, 
and $d_c=0.2$.
The fault is subjected by a uniform background shear $\tau^0=27.5~\text{MPa}$ and normal $\sigma^0=44~\text{MPa}$ stress.
We nucleate the fault rupture over a square region of size $a^2$ by instantaneously applying a loading traction of $31~\text{MPa}$, which locally exceeds the  peak friction strength.

\subsection{Results}

\begin{figure}[h]
    \centering
    \includegraphics[width=\textwidth]{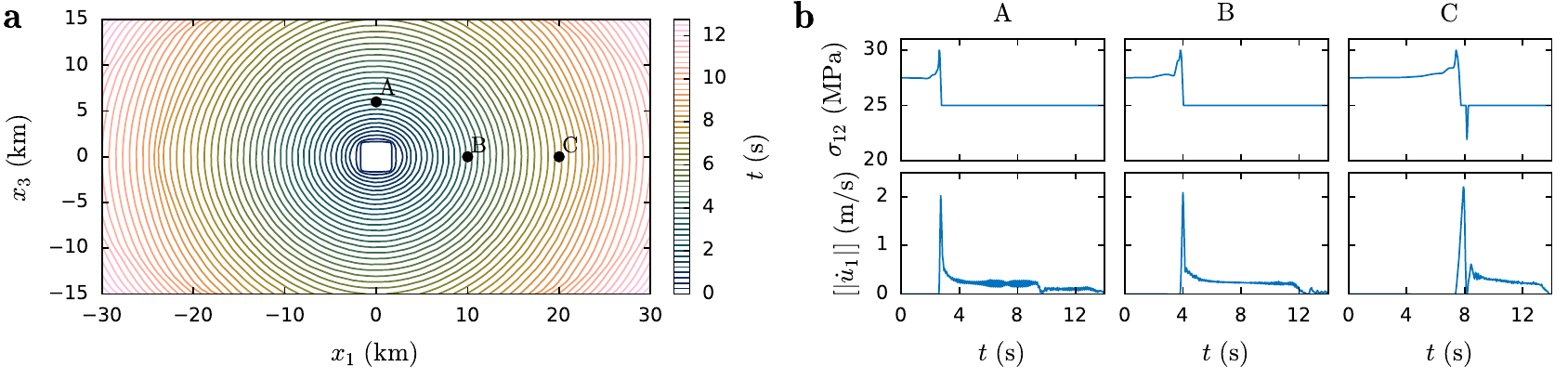}
    \caption{LVFZ setup solved using hybrid method with $\Delta x=100~\text{m}$.
    (a)~Contour of rupture front position each 0.5s.
    (b)~Fault shear stress, $\sigma_{12}$, and slip rate, $[|\dot{u}_1|]$ at three stations A, B, and C with position shown in (a). At station C the rupture has split into a slip-pulse followed by a crack-like rupture.}
    \label{fig:LVFZ_results}
\end{figure}

\begin{figure}[h]
    \centering
    \includegraphics[width=\textwidth]{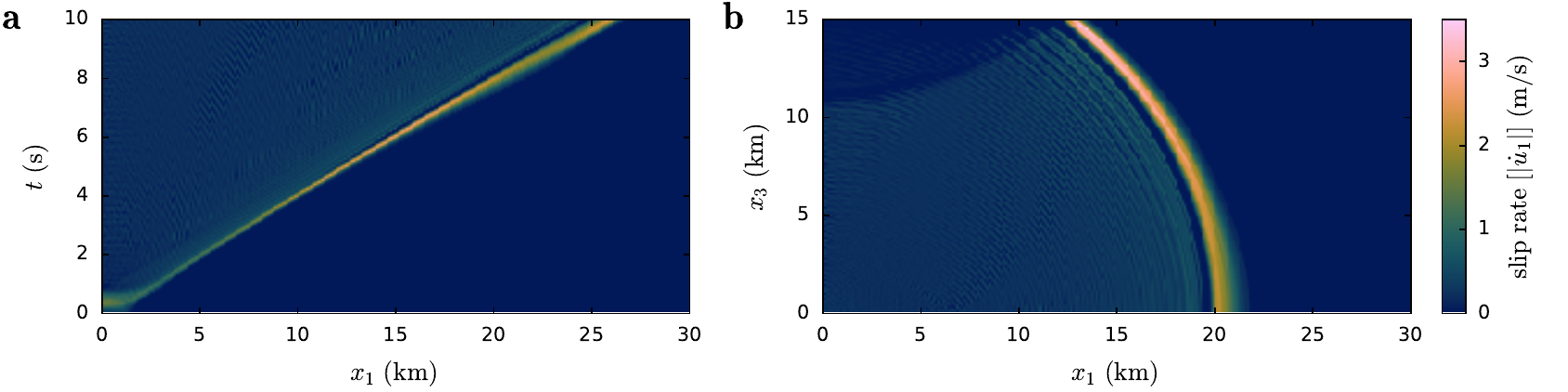}
    \caption{
      Pulse like rupture induced by LVFZ.
      (a)~Space-time diagram of slip rate $[|{\dot{u}_1}|]$ along the symmetry axis $x_3=0~\text{km}$. 
    At  $t\approx 5~\text{s}$, the rupture splits into a slip-pulse, followed by a crack-like rupture.
      (b)~Slip rate $[|{\dot{u}_1}|]$ at $t=8~\text{s}$ shows the spatial extent of the pulse-like and crack-like rupture.}
    \label{fig:LVFZ_results_pulse}
\end{figure}

As a result of the nucleation procedure, the rupture front quickly propagates radially and eventually spans the entire fault (see Fig.~\ref{fig:LVFZ_results}).
When a dynamic rupture propagates, it radiates elastic waves, which are then reflected at the boundary of the LVFZ.\cite{albertini_off-fault_2017}
Depending on the incident angle, the reflected wave can have an inverted polarization and cause unloading of the fault and generate a slip pulse.
This effect is also observed in our 3D simulations and is shown in Fig.~\ref{fig:LVFZ_results}b, station C, and in Fig.~\ref{fig:LVFZ_results_pulse}. The reflected wave causes the rupture to split into a pulse-like rupture, followed by a crack-like rupture.

The rapid acceleration and deceleration of a slip pulse are a source of high frequencies (see Fig.~\ref{fig:HVFZ_results}b Station C) and cause oscillations in slip velocity, trailing the rupture front (see Fig.~\ref{fig:LVFZ_results_pulse}). These oscillations do not affect the rupture propagation and disappear with further mesh refinement and regularized friction laws.\cite{kammer_existence_2014} Additionally, numerical damping, which is not used here, is often applied to minimize such high frequencies.
Since this problem cannot be solved with the SBI method, we validate the results of the hybrid method by varying the width of the virtual strip, $L_2$, and confirm that the solution is independent on the location of the elastodynamic boundary condition.

\section{Earthquake rupture in a heterogeneous medium: supershear transition}\label{sec:HVFZ}

\subsection{Setup}

The second showcase example, presented in this section, is similar to the previous example but with the more compliant material being the one at a distance from the fault. We consider a slip-weakening fault with an off-fault low velocity zone. This case could occur when a fault ruptures and interacts with the LVFZ of a nearby mature fault.
\cite{ma_effect_2015,albertini_off-fault_2017}
We use the same geometry, friction and elastic properties, and nucleation procedure as in Section~\ref{sec:LVFZ} but consider the fault to be embedded in the reference material, while beyond a distance $L^{FZ}/2$ from the fault plane the material has a {20\%} velocity reduction.

\subsection{Results}

\begin{figure}[h]
    \centering
    \includegraphics[width=\textwidth]{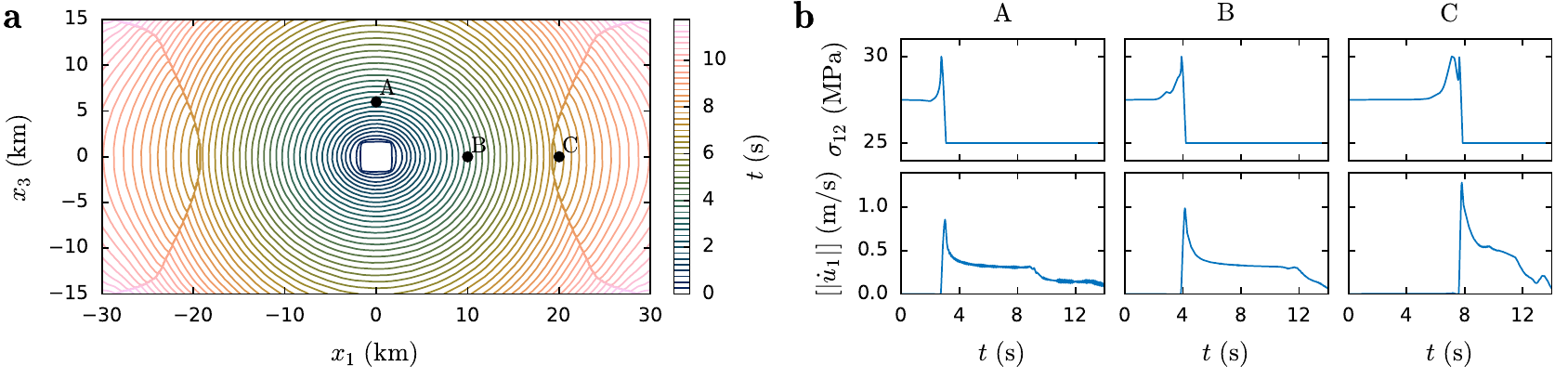}
    \caption{Off-fault low velocity zone setup solved using hybrid method with $\Delta x=100~\text{m}$. 
    (a)~Contour of rupture front position each 0.5s. Widely spaced contour lines represent supershear propagation region.
    (b)~Fault shear stress, $\sigma_{12}$, and slip rate, $[|\dot{u}_1|]$ at three stations A, B, and C with position shown in (a).}
    \label{fig:HVFZ_results}
\end{figure}

\begin{figure}[h]
    \centering
    \includegraphics[width=\textwidth]{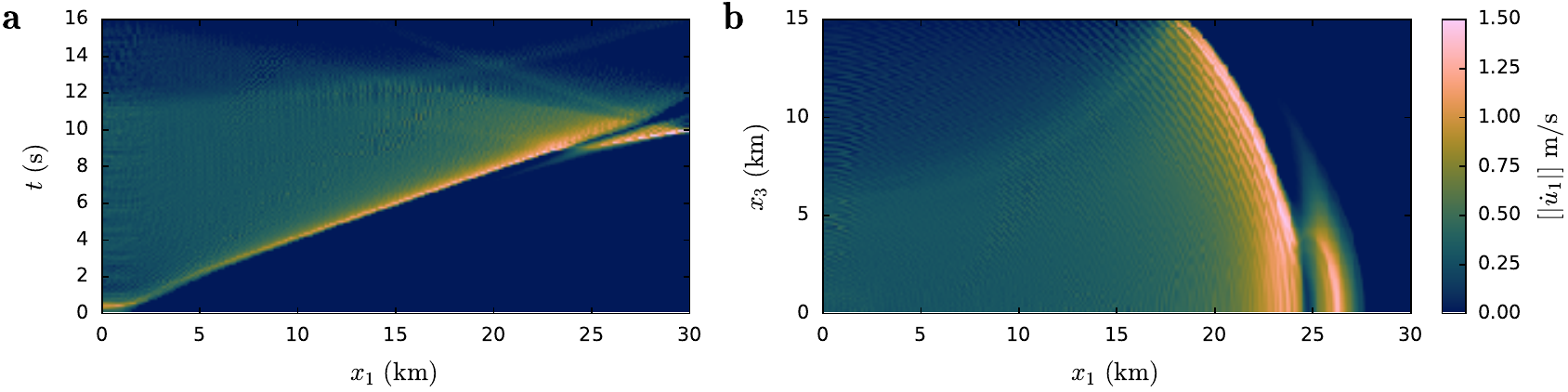}
    \caption{
      Supershear transition induced by off-fault low velocity zone.
      (a)~Space-time diagram of slip rate $[|\dot{u}_1|]$ along the symmetry axis $x_3=0~\text{km}$. The rupture transitions to supershear at $x_1\approx24~\text{km}$.
      (b)~Slip rate at $t=9.25~\text{s}$ shows the spatial extent of the supershear rupture, just after the transition has initiated.}
    \label{fig:HVFZ_results_super}
\end{figure}

Similar as in the LVFZ setup of Sec.~\ref{sec:LVFZ}, elastic waves are reflected at the boundary of the low velocity inclusion and affect the shear stress at the interface. 
However, the reflected waves have the same polarity as the incident ones and, hence, increase the shear stress in front of the propagating rupture front (see Fig.~\ref{fig:HVFZ_results}b). This increasing shear stress peak eventually causes the rupture to transition from subRayleigh to supershear velocities (see Fig.~\ref{fig:HVFZ_results}a and Fig.~\ref{fig:HVFZ_results_super}).
SubRayleigh propagation occurs when the rupture speed is lower than the Rayleigh wave speed, $c_R\approx 0.9 c_s$, and can be observed at stations A and B in Fig.~\ref{fig:HVFZ_results}.
Supershear propagation, however, refers to ruptures propagating faster than $c_s$ and their speed can approach the limiting speed, $c_p$.\cite{freund_mechanics_1979,kammer_equation_2018}
In our simulations, supershear rupture occurs within the domain surrounding station C in Fig.~\ref{fig:HVFZ_results}a.

In this 3D simulation, we observe a supershear transition through the Burridge-Andrews mechanism,\cite{burridge_admissible_1973,andrews_rupture_1976} where a shear stress peak in front of the existing crack nucleates the supershear rupture (see Fig.~\ref{fig:HVFZ_results}b, station B and   Fig.~\ref{fig:HVFZ_results_super}).
In contrast to 2D setups,\cite{ma_effect_2015,albertini_off-fault_2017} the extent of the supershear rupture is confined to a triangular shaped region, which surrounds station C. Additionally,  the transition occurs progressively: first at $x_3=0~\text{km}$ at $t\approx 7~\text{s}$, then it expands towards the $\pm e_3$ direction, and finally, at $t\approx 12~\text{s}$ it spans the entire seismogenic depth.

This example illustrates the ability of the hybrid method of successfully truncating the shear Mach front, radiated from the supershear rupture, without artificial reflections. Which allows us place the virtual boundary $S^\pm$ at only $2\Delta x$ from the the boundary of the low velocity inclusion.
As in the previous problem, we validate the results of the hybrid method by varying the width of the virtual strip, $L_2$. The solution is found to be independent on $L_2$.

\section{Earthquake rupture with step-over faults}\label{sec:step-over}

\subsection{Setup}

\begin{figure}[h]
    \centering
    \includegraphics[width=\textwidth]{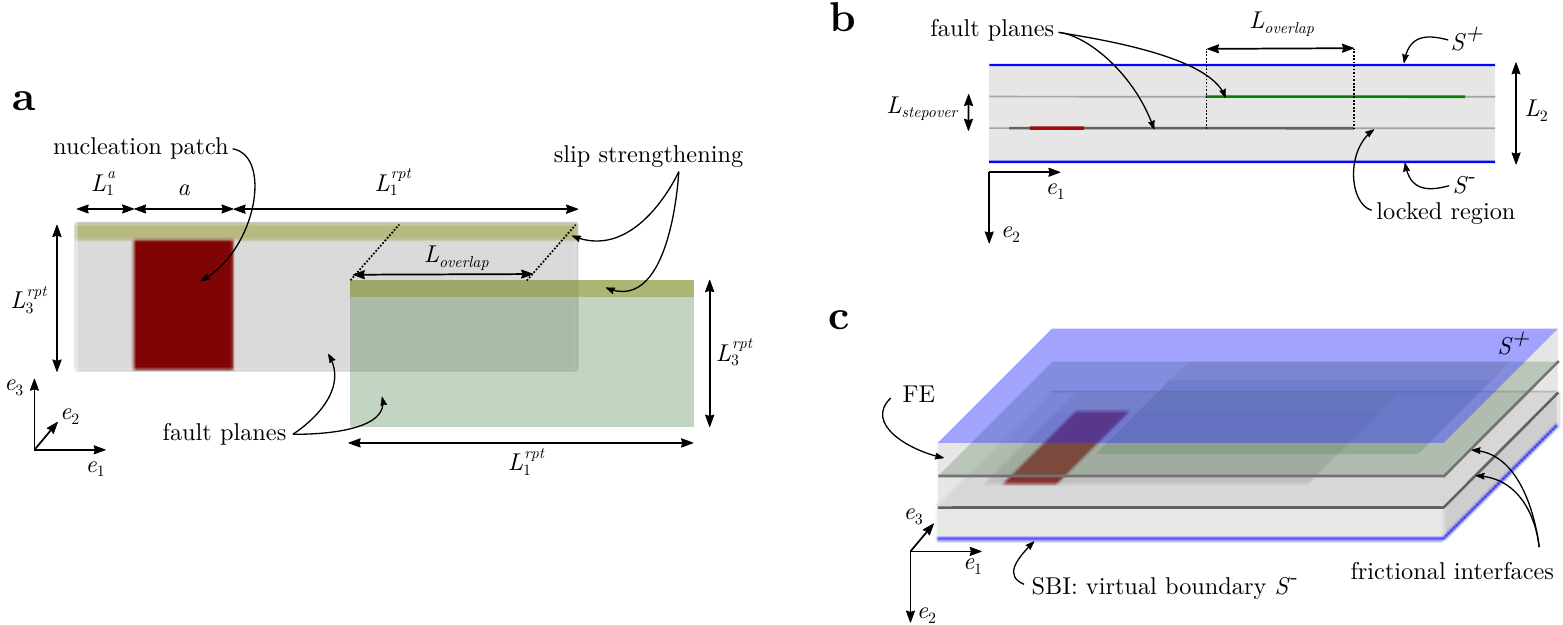}
    \caption{Setup of example earthquake rupture in unbounded domain with interacting parallel faults with step-over geometry.
    (a)~Setup geometry. Fault zone regions of rectangular size with $L_1^{rpt}=40~\text{km}$ and $L_3^{rpt}=10~\text{km}$, with nucleation over a width $a=20~\text{km}$ and the entire seismogenic depth, $L_3^{rpt}$. To the left of the nucleation patch the extent of  the primary fault is $L_1^a=10~\text{km}$.
    The two fault step-over geometry is characterized by  $L_{overlap}=20~\text{km}$ and $L_{step over}=1~\text{km}$, shown in (b).
    The faults are embedded in a homogeneous elastic medium.
    (b and c)~Hybrid setup: FE domain with SBI as elastodynamic boundary condition. The virtual strip width is $L_2=1.4~\text{km}$.
    }\label{fig:step-over_setup}
\end{figure}

Finally, we present an example of interaction between nearby faults, \textit{i.e.} two fracture planes side-by-side.
We consider a dilational step-over geometry with a system of two faults that overlap each other (see Fig.~\ref{fig:step-over_setup}). The dilational step-over implies that the location of the secondary fault with respect to the primary one is such that the rupture propagation on the primary fault will cause a temporary reduction in normal stress. 
The faults have uniform friction properties, $\mu_s=0.677$, $\mu_k=0.373$ and $d_c=0.5~\text{m}$, except on the top $1~\text{km}$ of the seismogenic zone, where a slip-strengthening condition is imposed.
At the bottom, \text{i.e.}, $x_3<L_3^{rpt}$, we consider a no slip boundary condition and the nucleation is achieved by an instantaneous reduction of the friction strength over a region of size $a\times L_3^{rpt}$ to it's kinetic value. This setup is analogous to a recent study by Bai and Ampuero.\cite{bai_effect_2017}
The elastic properties are the same as in Section~\ref{sec:TPV3}.
We choose a seismogenic depth, $L_3^{rpt}=10~\text{km}$ and uniform background shear $\tau_0=71.2~\text{MPa}$ and normal stress $\sigma_0=150~\text{MPa}$. Hence, the strength ratio is $S=(\mu_s \sigma_0 -\tau_0)/(\tau_0 - \mu_k \sigma_0)=1.75$, and the condition for the rupture to jump from one fault to the adjacent one is satisfied.\cite{bai_effect_2017}

Note that in this example there are no symmetries. In contrast to the previous examples where the fault plane represented a symmetry axis. As a result of the slip propagation on the primary fault, the normal stress on the secondary fault is not constant. The change in normal stress can reduce or increase the friction strength on the secondary fault, depending on their relative position. Thus, it can hinder or promote the rupture to jump between faults.

\subsection{Results}

\begin{figure}[h]
    \centering
    \includegraphics[width=\textwidth]{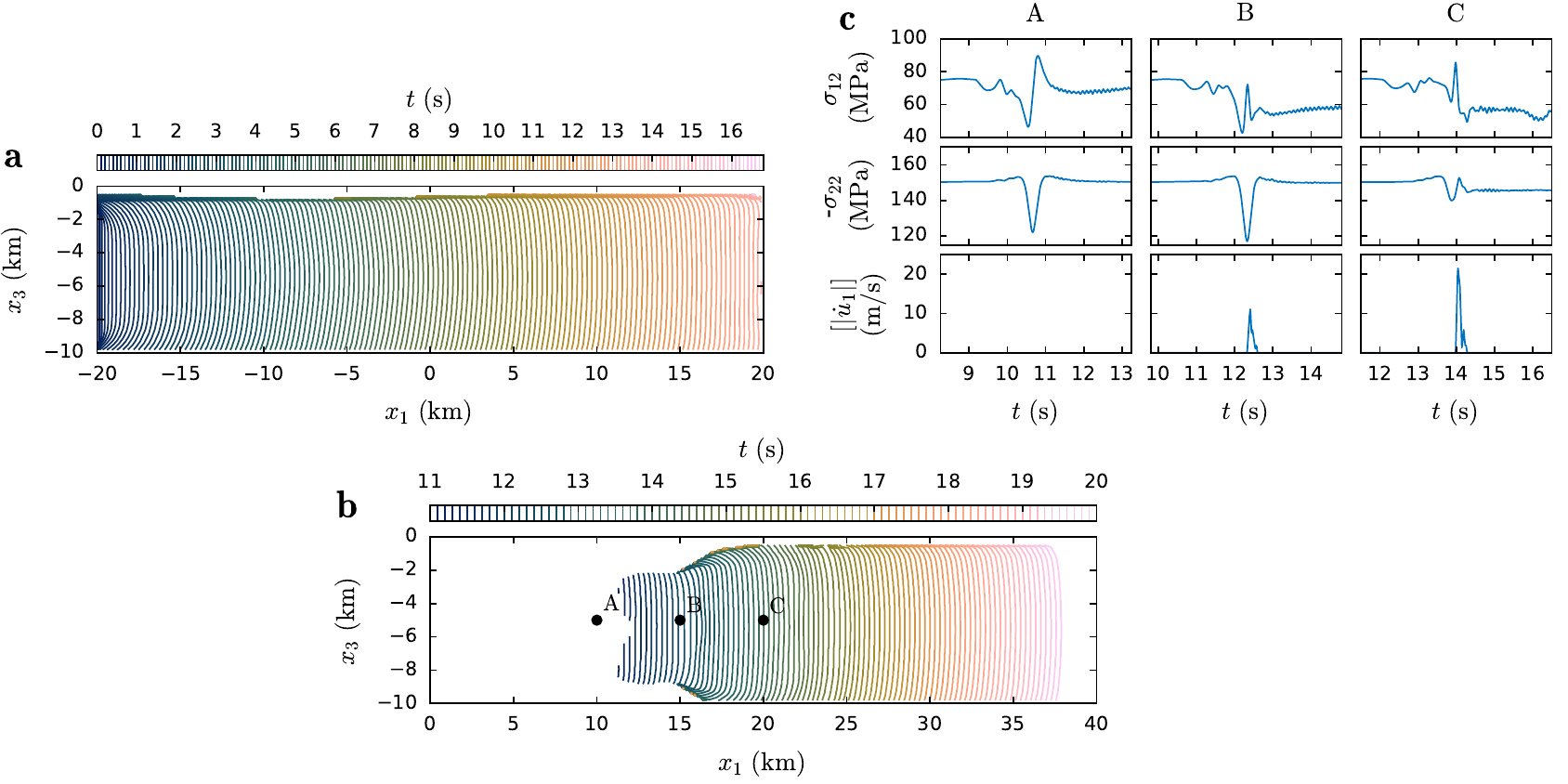}
    \caption{Results of example earthquake rupture in unbounded domain with interacting parallel faults with dilational step-over setup and spatial discretization $\Delta x=100~\text{m}$.
      (a)~Rupture front contour lines on primary fault, where nucleation occurs over region $-40~\text{km}<x_1<-20~\text{km}$ (not shown).
      (b)~Rupture front contour lines on secondary fault. Rupture transitions from primary to secondary fault.
      (c)~Fault shear, $\sigma_{12}$, and normal, $\sigma_{22}$, stress and slip rate, $[|\dot{u}_1|]$ at three stations A, B, and C  with position on the secondary fault shown in (b).}
      \label{fig:step-over_results}
\end{figure}

Using the hybrid method, we can successfully reproduce the results of Bai and Ampuero\cite{bai_effect_2017} (see Fig.~\ref{fig:step-over_results}):
after nucleation, the rupture propagates over the primary fault with a nearly vertical front, then the rupture jumps to the secondary fault and, as a  consequence of the no-slip boundary condition beyond the depth $L_3^{rpt}$, the rupture becomes a slip pulse.
The wave emitted by the primary rupture successfully nucleates a large rupture on the secondary fault in the forward direction.
This example illustrates the ability of the hybrid method to efficiently solve a large and complex simulation and efficiently truncate all incident waves without any artificial reflections.
To validate the results of the hybrid method we vary the width of the virtual strip, $L_2$, and find that also for this example $L_2$ does not affect the solution.

\section{Discussion}\label{sec:discussion}

We used a SCEC benchmark problem to validate the hybrid method and then demonstrated its flexibility and superior performance on more complex and heterogeneous problems.
The proposed hybrid method takes its flexibility to deal with nonlinearities or bulk heterogeneities from the FEM and its computational efficiency from the SBIM. In particular, since the SBIM provides a perfect wave absorption algorithm there is no artificial wave reflection at the virtual boundary.
Thus, one can reduce the width of the FEM domain arbitrarily close to the nonlinear or heterogeneous region, as long as the constitutive relation of the bulk beyond the virtual strip can be assumed to be linear.

\begin{figure}
    \centering
    \includegraphics[width=\textwidth]{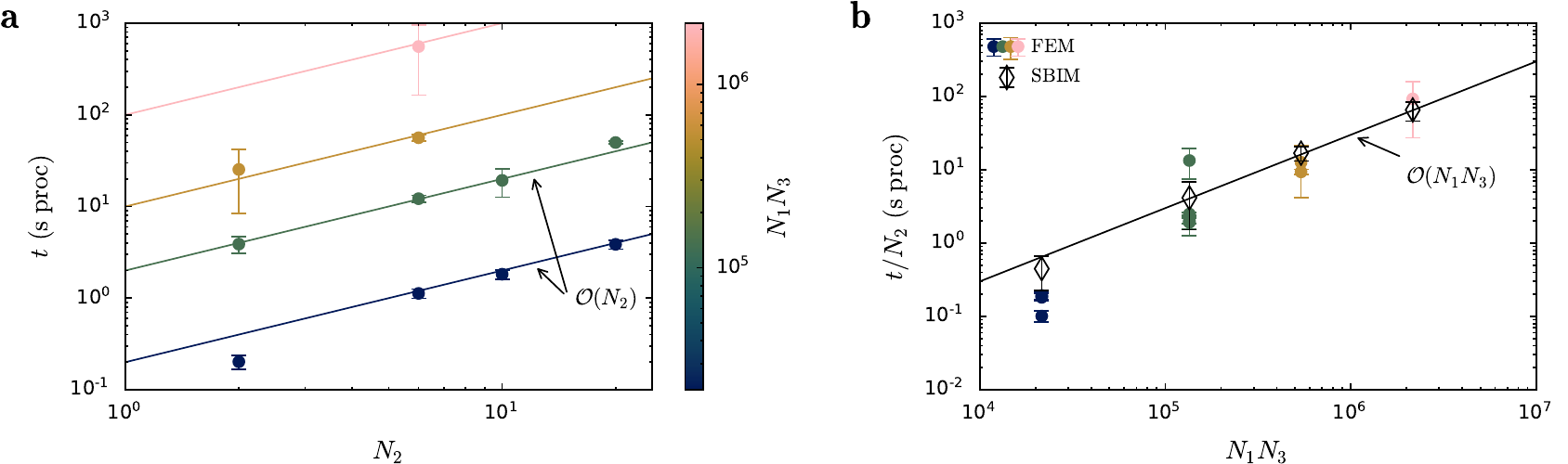}
    \caption{Performance study of FEM and Hybrid method. (a)~Computation time, $t$, of a FEM time step as function of width of the virtual strip, \textit{i.e.}, number of elements, $N_2$. Scaling of computational time is shown for a range of frictional interface discretizations, $N_1N_3$. 
    The complexity of the FEM time step is linear in $N_2$.
    (b)~Computation time of the FEM time step, $t$, normalized by $N_2$ (same data and color-code as in (a)). The computation time of a SBIM time step, when computing the elastodynamic boundary condition for the virtual strip, is equivalent to the computation time for one layer of FEM elements. Both, FEM and SBIM computation times are linear in $N_1N_3$.}
    \label{fig:performance}
\end{figure}

The computational savings of using the hybrid method -- instead of a traditional FEM -- can be assessed by considering the complexity of both FEM and SBIM.
The complexity of an explicit FEM time step is proportional to the number of degrees of freedom of the FEM problem,  \emph{i.e.,} $\mathcal{O}(N_1 N_2 N_3)$, where $N_i$ is the number of elements in the $i$-direction and we assume a regular mesh of hexagonal elements.
Similarly, the complexity of an SBIM time step scales with its number of degrees of freedom $\mathcal{O}(N_1 N_3)$. We measured the computation time for a range of simulations with different discretizations and domain sizes, which confirms the linear relationship between computational cost and the number of degrees of freedom (see Fig.~\ref{fig:performance}).
The computational saving of the hybrid method compared to a standard FEM lies in the reduction of $N_2$ due to the truncation of the FE domain.
Moreover, the added overhead cost of the SBI as wave absorption algorithm is in the same order of magnitude of only one layer of FE elements (see Fig.~\ref{fig:performance}b). Therefore, it is practically negligible.

For example, for a full FE simulation, $L_2$ must be in the order of $L_1$ to prevent artificial reflections at the domain boundary. However, using the hybrid method the domain size can be truncated up to the extent of the nonlinear region or the extent of the elastic heterogeneity, which are usually one to two orders of magnitude smaller than $L_1$.\cite{ma_effect_2015}
Assuming a regular spatial discretization, the domain truncation results in a reduction of $N_2$ by one order of magnitude and so will the computational cost. The savings may be even higher in other applications.

All our simulations were performed using distributed memory parallel computing with 48 threads and for the largest simulations 96 threads. 
Therefore, in Fig.~\ref{fig:performance}, we report the computational time multiplied by the number of parallel processes.
However, the SBIM library that we are using also supports shared memory parallelism and it is designed to be easily coupled to any FEM library written in C++. 
The only requirement is that the FEM mesh at the virtual boundary $S^\pm$ is a regular grid, due to the spectral representation of the boundary integral equations.

These computational savings represent an important step towards feasible modelling of complex temporal and spatial multi-scale 3D problems such as earthquake cycle simulations with near field heterogeneities, nonlinear material behavior and plasticity, as well as a networks of interacting faults, including fault branches and non-planar fault geometry.
The major challenge of earthquake cycle simulations is that they involve very long interseismic loading time (years) while the dynamic rupture happens extremely rapidly (seconds).
An advantage of the hybrid method is that the SBIM is already capable of absorbing elastic waves in the dynamic as well as in the quasi-dynamic limit and these approaches can be combined in a variable time stepping scheme, introduced by Lapusta et al.\cite{lapusta_elastodynamic_2000}
Such a temporal multi-scale simulation couples a quasi-dynamic SBIM with an implicit FEM during the slow loading phase and, once the ruptures become dynamic, it switches to a dynamic SBIM coupled with an explicit FEM -- as considered in the current study.
These variable time-stepping hybrid method was introduced in a 2D antiplane framework\cite{abdelmeguid_novel_2019} and will be extended to 3D in future work.

Another advantage of the hybrid method is that it could be implemented with any volume based method. For example, if the fault plane is not known a priori, one could use discretization techniques with embedded discontinuities, such as the XFEM\cite{liu_extended_2009} or the discontinuous Galerkin.\cite{pelties_three-dimensional_2012}

\section{Conclusion}\label{sec:conclusion}

We developed a three dimensional hybrid method  combining the finite-element method with the spectral boundary-integral method.
We validated the hybrid method using a benchmark problem and illustrated its potential for solving complex earthquake propagation on various example problems including near systems with field heterogeneity and  multiple interacting faults.
The hybrid method is suitable for cases where the spatial extent of near field nonlinearity and heterogeneity is too large to be lumped into an effective fault constitutive law, but is still considerably smaller than the domain of interest for the wave propagation.
In these cases, the hybrid method allows for a reduction of computational cost by at least one order of magnitude with respect to a full finite-element implementation, while maintaining the same level of accuracy.
The high accuracy and computational efficiency of the hybrid method enable the investigation of complex failure problems such as multi-physics fault zone problems.

\section*{Acknowledgments}
The authors thank Dr. Nicolas Richart for technical help. 


\end{document}